\documentclass{article} % use \documentstyle for old LaTeX compilers

\usepackage[utf8]{inputenc} % 'cp1252'-Western, 'cp1251'-Cyrillic, etc.
\usepackage[english]{babel} % 'french', 'german', 'spanish', 'danish', etc.
\usepackage{amsmath}
\usepackage{amssymb}
\usepackage{txfonts}
\usepackage{mathdots}
\usepackage[classicReIm]{kpfonts}
\usepackage{graphicx}

\title{Fractional integration and differentiation}
\author{Yachmenev Andrey Igorevich\\ Moscow State University of Technology "Stankin" \\ Vadkovsky Lane, 1, Moscow, 127055, Russia,\\  Undergraduate student\and Yaremko Oleg Emanuilovich\\ Moscow State University of Technology "Stankin"\\  Vadkovsky Lane, 1, Moscow,  127055, Russia,\\  Doctor of Physics and Mathematics, Professor, Docent}

\begin{document}
\maketitle
\section*{Abstract}
%\selectlanguage{english} % remove comment delimiter ('%') and select language if required

In this paper, we introduce a new method for calculating fractional integrals and differentials. The method involves an equation that we have obtained from infinite applied integration by parts. The equation works for special class of functions and provides a series representation of integration.  This representation will be useful for working with smooth functions and for approximation due to the potential reduction of the value of the sequential elements.

\noindent Keywords: Fractional integral, fractional differential, approximation, differential equation.

\section{Introduction}

\noindent 

Fractional integrals and differentials were first studied in the early 19th century by mathematicians Liouville and Riemann. They were later expanded upon by several other scientists, including Cauchy, Laplace and Lebesgue [1]. The use of fractional differentials is now seen in various fields such as diffusion, thermal conductivity, fractal structures, and much more [1]. In modern representation of fractional calculus, there are multiple definitions, but two main ones stand out [2].

Riemann--Liouville integral:
\begin{equation} \label{1.1} 
I^{\alpha } f(x)=\frac{1}{\Gamma (\alpha )} \int _{\alpha }^{x}f(t)(x-t)^{\alpha -1} dt  
\end{equation} 
\begin{equation} \label{1.2} 
D_{x}^{\alpha } f(x)=\left\{\begin{array}{c} {\frac{d^{\left\lceil \alpha \right\rceil } }{dx^{\left\lceil \alpha \right\rceil } } I_{x}^{\left\lceil \alpha \right\rceil -\alpha } f(x)\quad \alpha >0} \\ {f(x)\quad \alpha =0} \\ {I_{x}^{-\alpha } f(x)\quad \alpha <0} \end{array}\right\} 
\end{equation}

Grunwald-Letnikov derivative:

left-sided derivative:
\begin{equation} \label{1.3} 
{}^{GL} D_{a^{+} }^{\alpha } f(x)={\mathop{\lim }\limits_{h\to 0}} \frac{1}{h^{\alpha } } \sum _{k=0}^{\left\lfloor n\right\rfloor }(-1)^{k} \frac{\Gamma (\alpha +1)f(x-kh)}{\Gamma (k+1)\Gamma (\alpha -k+1)}  ,\quad nh=x-a 
\end{equation}

right-sided derivative:
\begin{equation} \label{1.4} 
{}^{GL} D_{b^{-} }^{\alpha } f(x)={\mathop{\lim }\limits_{h\to 0}} \frac{1}{h^{\alpha } } \sum _{k=0}^{\left\lfloor n\right\rfloor }(-1)^{k} \frac{\Gamma (\alpha +1)f(x+kh)}{\Gamma (k+1)\Gamma (\alpha -k+1)}  ,\quad nh=b-x 
\end{equation} 

Other definitions attempt to solve certain specific cases, but they are conceptually similar to either the first or the second. The peculiarity of calculating these approaches is similar in that it requires to calculate the sum of all elements in some approximation, since each of them can have an equal influence on the result.

\section{Derivation of a new equation}

The basis for the derivation of a new representation of the fractional derivative will be one of the basic properties of integration - integration by parts. For our objectives, we define the boundaries of integration from the constant (zero) to the variable x [3].
\begin{equation} \label{1.5} 
\int _{0}^{x}f(t)dt=\left|\begin{array}{cc} {f(t)=u} & {dt=dv} \\ {f'(t)=du} & {t=v} \end{array}\right| =f(x)x-\int _{0}^{x}f'(t)tdt  
\end{equation} 

After applying this method, we obtain two components: an expression that does not contain an integral, and an expression similar to the original integral. Next, we will use integration by parts sequentially on the integral that appears at each step.
\begin{equation} \label{1.6} 
\begin{array}{l} {\int _{0}^{x}f(t)dt =f(x)x-\int _{0}^{x}f'(t)tdt= f(x)x-\frac{f'(x)x^{2} }{2} +\int _{0}^{x}\frac{f''(t)t^{2} }{2}  dt=\ldots } \\ {=\sum _{k=0}^{n}\frac{(-1)^{k} f^{(k)} (x)x{}^{k+1} }{(k+1)!}  +(-1)^{n+1} \int _{0}^{x}\frac{f^{(n+1)} (t)t{}^{n+1} }{(n+1)!}  dt} \end{array} 
\end{equation} 

After that, we turn to the limiting transition of the two resulting parts with $n\to \infty $. To obtain a new expression, it is worth highlighting a class of functions for which the integral will tend to zero at the limit transition.

A suitable class of functions can be defined as$f(x)\in C^{\infty } \vee {\mathop{\lim }\limits_{n\to \infty }} \frac{f^{(n)} (x)}{n!} =0$

For this class of functions, the integral can be represented as an infinite sum of their derivatives.
\begin{equation} \label{1.7} 
\int _{0}^{x}f(t)dt= {\mathop{\lim }\limits_{n\to \infty }} (\sum _{k=0}^{n}\frac{(-1)^{k} f^{(k)} (x)x{}^{k+1} }{(k+1)!}  +(-1)^{n+1} \int _{0}^{x}\frac{f^{(n+1)} (t)t{}^{n+1} }{(n+1)!}  dt)=\sum _{k=0}^{\infty }\frac{(-1)^{k} f^{(k)} (x)x{}^{k+1} }{(k+1)!}   
\end{equation} 

To utilize a different constant for integration, it is necessary to make a slight adjustment to the equation.
\begin{equation} \label{1.8} 
\begin{array}{l} {\int _{\alpha }^{x}f(x)dx =\int _{0}^{x}f(x)dx-\int _{0}^{\alpha }f(x)dx} \\  {=\sum _{k=0}^{\infty }\frac{(-1)^{k} }{(k+1)!} (x^{k+1} \frac{d^{k} f}{dx^{k} }  (x)-\alpha ^{k+1} \frac{d^{k} f}{dx^{k} } (a)) } \end{array}
\end{equation} 

After getting Eq. \eqref{1.7} and Eq. \eqref{1.8} we present some illustrative examples.

\noindent \textbf{Example 1.}

Finding power rule for integration. 
\[\int _{0}^{x}t^{m} dt=\frac{x^{m+1} }{m+1}  \] 
\[\sum _{k=0}^{\infty }\frac{(-1)^{k} x{}^{k+1} }{(k+1)!} \frac{d^{k} }{dx^{k} }  x^{m} =\sum _{k=0}^{m}\frac{(-1)^{k} x^{k+1} m!x^{m-k} }{(k+1)!(m-k)!}  =x^{m+1} \sum _{k=0}^{m}\frac{(-1)^{k} m!}{(k+1)!(m-k)!}  =\frac{x^{m+1} }{m+1} \] 
\textbf{Example 2.}

Finding integration of exponent.
\[\int _{0}^{x}e^{t} dt= e^{t} -1\] 
\[\sum _{k=0}^{\infty }\frac{(-1)^{k} x{}^{k+1} }{(k+1)!} \frac{d^{k} }{dx^{k} }  e^{x} =\sum _{k=0}^{\infty }\frac{(-1)^{k} x{}^{k+1} }{(k+1)!}  e^{x} =e^{x} \sum _{k=0}^{\infty }\frac{(-1)^{k} x{}^{k+1} }{(k+1)!}  =e^{x} (-e^{-x} +1)=e^{x} -1\] 
\textbf{Example 3.}

Finding integration of natural logarithm.
\[\int _{1}^{x}\ln (x)dx=x\ln (x)-x+1 \] 
\[\sum _{k=0}^{\infty }\frac{(-1)^{k} }{(k+1)!} (x^{k+1} \frac{d^{k} \ln }{dx^{k} }  (x)-\alpha ^{k+1} \frac{d^{k} \ln }{dx^{k} } (\alpha ))=x\ln (x)-0+\sum _{k=0}^{\infty }\frac{(-1)^{k+1} }{(k+2)!} (x^{k+2} \frac{d^{k} x^{-1} }{dx^{k} }  -1^{k+2} \frac{d^{k} x^{-1} }{dx^{k} } (1))\] 
\[=x\ln (x)+\sum _{k=0}^{\infty }\frac{(-1)^{k+1} }{(k+2)!} (x^{k+2} (-1)^{k} k!x^{-1-k}  -(-1)^{k} k!)=x\ln (x)-\sum _{k=0}^{\infty }\frac{1}{(k+1)(k+2)} (x -1)\] 
\[=x\ln (x)-\sum _{k=0}^{\infty }(\frac{1}{k+1} -\frac{1}{k+2} )(x -1)=x\ln (x)-x+1\] 
\section{Generalization}

The equation for double integration can be derived from the Eq. 1.7.

\noindent  Let's suppose ${\rm F(x)\; =\; }\int _{0}^{x}f(x)dx $
\begin{equation} \label{1.9} 
\begin{array}{l} {\int _{0}^{x}\int _{0}^{x}f(x)dxdx=\int _{0}^{x}F(x)dx= \sum _{k=0}^{\infty }\frac{(-1)^{k} x{}^{k+1} }{(k+1)!} \frac{d^{k} }{dx^{k} }    F(x)=} \\ {=x\int _{0}^{x}f(x)dx +\sum _{k=1}^{\infty }\frac{(-1)^{k} x{}^{k+1} }{(k+1)!} \frac{d^{k-1} }{dx^{k-1} }  f(x)} \\ {=\sum _{k=0}^{\infty }\frac{(-1)^{k} x{}^{k+2} }{(k+1)!} \frac{d^{k} }{dx^{k} }  f(x)+\sum _{k=0}^{\infty }\frac{(-1)^{k+1} x{}^{k+2} }{(k+2)!} \frac{d^{k} }{dx^{k} }  f(x)} \\ {=\sum _{k=0}^{\infty }(\frac{(-1)^{k} }{(k+1)!} (1+\frac{(-1)}{(k+2)} )x^{k+2} \frac{d^{k} }{dx^{k} }  f(x)=} \\ {=\sum _{k=0}^{\infty }(\frac{(-1)^{k} (k+1)}{(k+2)!} x^{k+2} \frac{d^{k} }{dx^{k} }  f(x)} \end{array} 
\end{equation} 

Using the same technique, it is possible to obtain an equation for a larger number of integrals applied in a row. In general, it is possible to derive a new equation for n consecutive integrals applied.
\begin{equation} \label{1.10} 
f^{(-n)} (x)=\frac{1}{(n-1)!} \sum _{k=0}^{\infty }\frac{(-1)^{k} x^{k+n} }{k!(k+n)} \frac{d^{k} }{dx^{k} } f(x)  
\end{equation} 
Eq. \eqref{1.10} can be proved by induction.

\noindent 1.Base case we already get.
\[f^{(-1)} (x)=\sum _{k=0}^{\infty }\frac{(-1)^{k} x^{k+1} }{k!(k+1)} \frac{d^{k} }{dx^{k} } f(x) \] 
2. Induction step

\noindent Let's suppose ${\rm F(x)\; =\; }\int _{0}^{x}f(x)dx $
\[\begin{array}{l} {f^{(-n-1)} (x)=F^{(-n)} (x)=\frac{1}{(n-1)!} \sum _{k=0}^{\infty }\frac{(-1)^{k} x^{k+n} }{k!(k+n)} \frac{d^{k} }{dx^{k} } F(x)=\frac{F(x)x^{n} }{n!} +\frac{1}{(n-1)!} \sum _{k=1}^{\infty }\frac{(-1)^{k} x^{k+n} }{k!(k+n)} \frac{d^{k} }{dx^{k} } F(x)  =} \\ {\frac{x^{n} }{n!} \sum _{k=0}^{\infty }\frac{(-1)^{k} x{}^{k+1} }{(k+1)!} \frac{d^{k} }{dx^{k} }  f(x)+\frac{1}{(n-1)!} \sum _{k=1}^{\infty }\frac{(-1)^{k} x^{k+n} }{k!(k+n)} \frac{d^{k} }{dx^{k} } F(x)= \sum _{k=0}^{\infty }\frac{(-1)^{k} x{}^{k+n+1} }{n!(k+1)!} \frac{d^{k} }{dx^{k} }  f(x)+} \\ {\frac{1}{(n-1)!} \sum _{k=0}^{\infty }\frac{(-1)^{k+1} x^{k+n+1} }{(k+1)!(k+n+1)} \frac{d^{k} }{dx^{k} } f(x)= \sum _{k=0}^{\infty }\frac{(-1)^{k} x{}^{k+n+1} }{(k+1)!(n-1)!} (\frac{1}{n} +\frac{-1}{(n+k+1)} )\frac{d^{k} }{dx^{k} }  f(x)=} \\ {\sum _{k=0}^{\infty }\frac{(-1)^{k} x{}^{k+n+1} }{(k+1)!(n-1)!} \frac{k+1}{n(n+k+1)} \frac{d^{k} }{dx^{k} }  f(x)=\frac{1}{((n+1)-1)!} \sum _{k=0}^{\infty }\frac{(-1)^{k} x{}^{k+(n+1)} }{(k+1)!(k+(n+1))} \frac{d^{k} }{dx^{k} }  f(x)=f^{(-n-1)} (x)} \end{array}\] 
Q.E.D.

Now, we can go into the space of fractional integration, since all the elements adhere to the analytic continuation, and derive an expression of the following form:
\begin{equation} \label{1.11} 
{}_{0} I_{x}^{\alpha } f(x)=\frac{1}{\Gamma (\alpha )} \sum _{k=0}^{\infty }\frac{(-1)^{k} x^{k+\alpha } }{k!(k+\alpha )} \frac{d^{k} }{dx^{k} } f(x)  
\end{equation} 

Further, using the definition \eqref{1.2}, we can obtain the fractional differential equation.

\noindent Let's suppose$n=\left\lceil \alpha \right\rceil $ for $\alphaup$ $\mathrm{>}$ 0
\begin{equation} \label{1.12} 
\begin{array}{l} {{}_{0} D_{x}^{\alpha } f(x)=\frac{d^{n} }{dx^{n} } {}_{0} I_{x}^{n-\alpha } f(x)=\frac{d^{n} }{dx^{n} } \frac{1}{\Gamma (n-\alpha )} \sum _{k=0}^{\infty }\frac{(-1)^{k} x^{k+n-\alpha } }{k!(k+n-\alpha )} \frac{d^{k} }{dx^{k} } f(x)= } \\ {=\frac{1}{\Gamma (n-\alpha )} \sum _{k=0}^{\infty }\frac{(-1)^{k} }{k!(k+n-\alpha )} \frac{d^{n} }{dx^{n} } (x^{k+n-\alpha } \frac{d^{k} }{dx^{k} } f(x))= } \\ {=\frac{1}{\Gamma (n-\alpha )} \sum _{k=0}^{\infty }\frac{(-1)^{k} }{k!(k+n-\alpha )} \sum _{b=0}^{n}\left(\begin{array}{c} {n} \\ {b} \end{array}\right) \frac{d^{n-b} }{dx^{n-b} } (x^{k+n-\alpha } )\frac{d^{k+b} }{dx^{k+b} } f(x) =} \\ {=\frac{1}{\Gamma (n-\alpha )} \sum _{k=0}^{\infty }\frac{(-1)^{k} }{k!(k+n-\alpha )} \sum _{b=0}^{n}\left(\begin{array}{c} {n} \\ {b} \end{array}\right) \frac{\Gamma (k+n-\alpha +1)}{\Gamma (k+b-\alpha +1)} x^{k+b-\alpha } \frac{d^{k+b} }{dx^{k+b} } f(x) =} \\ {=\frac{1}{\Gamma (n-\alpha )} \sum _{b=0}^{n}\left(\begin{array}{c} {n} \\ {b} \end{array}\right) \sum _{k=0}^{\infty }\frac{(-1)^{k} }{k!} \frac{\Gamma (k+n-\alpha )}{\Gamma (k+b-\alpha +1)} x^{k+b-\alpha } \frac{d^{k+b} }{dx^{k+b} } f(x) =} \\ {=\sum _{b=0}^{n}\left(\begin{array}{c} {n} \\ {b} \end{array}\right) \sum _{k=0}^{\infty }(-1)^{k} \left(\begin{array}{c} {k+n-\alpha -1} \\ {k} \end{array}\right)\frac{1}{\Gamma (k+b-\alpha +1)} x^{k+b-\alpha } \frac{d^{k+b} }{dx^{k+b} } f(x) =} \\ {=\sum _{b=0}^{n}\left(\begin{array}{c} {n} \\ {b} \end{array}\right) \sum _{k=0}^{\infty }\left(\begin{array}{c} {\alpha -n} \\ {k} \end{array}\right)\frac{x^{k+b-\alpha } }{\Gamma (k+b-\alpha +1)} \frac{d^{k+b} }{dx^{k+b} } f(x) } \end{array} 
\end{equation} 

 If we equalize equation \eqref{1.12} according to the degree of differentiation, we get an equation similar to equation \eqref{1.11}
\begin{equation} \label{1.13} 
{}_{0} D_{x}^{\alpha } f(x)=\sum _{k=0}^{\infty }\left(\begin{array}{c} {\alpha } \\ {k} \end{array}\right)\frac{x^{k-\alpha } }{\Gamma (k-\alpha +1)} \frac{d^{k} }{dx^{k} } f(x)  
\end{equation} 

This equation work for any complex $\mathrm{\alphaup }$.

\section{Operator scope}

It is a complex task to investigate the convergence of this series. The basis of convergence of the series should be expanded in comparison to the previous definition. Now it is worth defining a class of functions as $f(x)\in C^{\infty } \vee {\mathop{\lim }\limits_{n\to \infty }} \left(\begin{array}{c} {\alpha } \\ {n} \end{array}\right)\frac{f^{(n)} (x)}{\Gamma (n+1-\alpha )} =0$

It is also important to introduce a restriction on the growth of the function. Using Fabry's theorem and Field's theorem, we can limit the equation \eqref{1.13} to the sum with gaps [4][5].

${}_{\ }$Let's suppose $N_{n} =\{ \omega _{n} \in {\mathbb R}|\exists h\in {\mathbb R}:\left(\begin{array}{c} {\alpha } \\ {n} \end{array}\right)\frac{f^{(n)} (z)}{\Gamma (n+1-\alpha )} \le hz^{\omega _{n} } \} $

So, if for every n $N_{n} \ne \emptyset $ and sequence $\min \left(N_{n} \right)$ grows linearly, sum converge. 

The radius of convergence can be determined by standard methods after function substitution.

\section{Properties}

Now it is possible to define the properties of the new equation. To describe some properties, it is useful to define a falling factorial [6].
\begin{equation} \label{1.14} 
\left(x\right)^{\underline{n}} =\prod _{k=1}^{n}(x-(k-1))  
\end{equation}

Derivative of the sum:
\begin{equation} \label{1.15} 
\begin{array}{l} {{}_{0} D_{x}^{\alpha } [f\left(x\right)+g(x)]=\sum _{k=0}^{\infty }\left(\begin{array}{c} {\alpha } \\ {k} \end{array}\right)\frac{x^{k-\alpha } }{\Gamma (k-\alpha +1)} \frac{d^{k} }{dx^{k} } [f(x) +g(x)]} \\ {=\sum _{k=0}^{\infty }\left(\begin{array}{c} {\alpha } \\ {k} \end{array}\right)\frac{x^{k-\alpha } }{\Gamma (k-\alpha +1)} \frac{d^{k} }{dx^{k} } f(x) +\sum _{k=0}^{\infty }\left(\begin{array}{c} {\alpha } \\ {k} \end{array}\right)\frac{x^{k-\alpha } }{\Gamma (k-\alpha +1)} \frac{d^{k} }{dx^{k} } g(x) } \\ {={}_{0} D_{x}^{\alpha } f\left(x\right)+{}_{0} D_{x}^{\alpha } g(x)} \end{array} 
\end{equation}

Derivative of the product:
\begin{equation} \label{1.16} 
\begin{array}{l} {{}_{0} D_{x}^{\alpha } (f(x)g(x))=\sum _{k=0}^{\infty }\left(\begin{array}{c} {\alpha } \\ {k} \end{array}\right)\frac{x^{k-\alpha } }{\Gamma (k-\alpha +1)} \frac{d^{k} }{dx^{k} } [f(x) g(x)]} \\ {=\sum _{k=0}^{\infty }\left(\begin{array}{c} {\alpha } \\ {k} \end{array}\right)\frac{x^{k-\alpha } }{\Gamma (k-\alpha +1)} \sum _{h=0}^{k}\left(\begin{array}{c} {k} \\ {h} \end{array}\right) \frac{d^{h} }{dx^{h} } f(x)\frac{d^{k-h} }{dx^{k-h} } g(x)= } \\ {=\sum _{k=0}^{\infty }\frac{\left(\alpha \right)^{\underline{k}} x^{k+\alpha } }{\Gamma (k-\alpha +1)} \sum _{h=0}^{k}\frac{1}{h!}  \frac{d^{h} }{dx^{h} } f(x)\frac{1}{(k-h)!} \frac{d^{k-h} }{dx^{k-h} } g(x) =} \\ {=\frac{\Gamma (\alpha +1)}{\pi } \sum _{k=0}^{\infty }\frac{\sin (\pi \alpha -\pi k)x^{k+\alpha } }{(\alpha -k)} \sum _{h=0}^{k}\frac{1}{h!}  \frac{d^{h} }{dx^{h} } f(x)\frac{1}{(k-h)!} \frac{d^{k-h} }{dx^{k-h} } g(x) } \end{array} 
\end{equation}

The sum of differentials of different degrees of one function
\begin{equation} \label{1.17} 
\begin{array}{l} {a{}_{0} D_{x}^{\alpha } f(x)+b{}_{0} D_{x}^{\beta } f(x)} \\ {=\sum _{k=0}^{\infty }\left(\begin{array}{c} {\alpha } \\ {k} \end{array}\right)\frac{x^{k-\alpha } }{\Gamma (k-\alpha +1)} \frac{d^{k} }{dx^{k} } f(x)+\sum _{k=0}^{\infty }\left(\begin{array}{c} {\beta } \\ {k} \end{array}\right)\frac{x^{k-\beta } }{\Gamma (k-\beta +1)} \frac{d^{k} }{dx^{k} } f(x)  =} \\ {=\sum _{k=0}^{\infty }\frac{x^{k} }{k!} (\frac{\left(\alpha \right)^{\underline{k}} x^{\alpha } }{\Gamma (k-\alpha +1)} +\frac{\left(\beta \right)^{\underline{k}} x^{\beta } }{\Gamma (k-\beta +1)} )\frac{d^{k} }{dx^{k} } f(x) } \end{array} 
\end{equation}
\section{Application}

The new equation possesses several distinguishing features that make it practical for application compared to equations \eqref{1.1}, \eqref{1.2} and \eqref{1.3}, \eqref{1.4}. Firstly, its structure is similar to the Taylor series, leading to a decrease of each subsequent element decreases inversely proportionally to the factorial during numerical calculation. Secondly, it does not require special algorithms for extending into a set of complex numbers. By utilizing equation \eqref{1.13}, it is possible to devise an algorithm for approximate calculation. The first version of the algorithm will be presented in the Wolfram Mathematica language.

\begin{tabular}{|p{4.4in}|} \hline 
$\mathrm{FractDiff}[\mathrm{pow\_},\mathrm{f\_},\mathrm{t\_},\mathrm{noo\_}]:=\mathrm{Sum}[\mathrm{Binomial}[\mathrm{pow},k]*{t\wedge (k-\mathrm{pow})}/{\mathrm{Gamma}[k+1-\mathrm{pow}]}*D[f,\{t,k\}],\{k,0,\mathrm{noo}\}]$ \\ \hline 
\end{tabular}

The function takes 4 arguments as input:

\noindent 1. pow (power) - Power of the differentiation operator.

\noindent 2. f (function) - The function to which the operator is applied.

\noindent 3. t - The variable by which calculations are performed.

\noindent 4. noo (Number of Operations) - The number of operations that determine the accuracy of calculations.

This version will be useful for scientific calculations. However, for use in large projects, a C++ implementation may be more practical. Equation \eqref{1.4} will be used to calculate the integer power of the derivative of function in C++. By combining these approaches, though with slightly increased time costs, the accuracy of calculations increases greatly [7].

We will need auxiliary functions that are not contained in standard libraries:
\begin{enumerate}
\item Binomiald - a function that calculates the binomial coefficient, capable of calculating non-integer values

\begin{tabular}{|p{4.4in}|} \hline 
double binomiald(double k, int a) noexcept\newline $\mathrm{\{}$\newline  double b = 1.0;\newline  double k2 = k + 1.0;\newline  double p = 1.0;\newline  for (int c = 1; c $\mathrm{<}$= a; ++c) $\mathrm{\{}$\newline   b *= k2 / p - 1.0;\newline   p += 1.0;\newline  $\mathrm{\}}$\newline  return b;\newline $\mathrm{\}}$ \\ \hline 
\end{tabular}

\item  Grunwald-Letnikov derivative

\begin{tabular}{|p{4.4in}|} \hline 
double GL\_righthand\_derivative(std::function$\mathrm{<}$double(double)$\mathrm{>}$ f, double pow, double point, unsigned long N) noexcept\newline $\mathrm{\{}$\newline     const double h = point / N;\newline     double sign = 1.0;\newline     double sum = 0;\newline     const double correction = pow * h / 2.0;\newline     double pointshift = point - correction - h;\newline     for (size\_t k = 0; k $\mathrm{<}$ N; k++)//gofor(k, range(0, N))\newline     $\mathrm{\{}$\newline         sum += sign * binomiald(pow, k) * f(pointshift += h);\newline         sign *= -1.0;\newline     $\mathrm{\}}$\newline     return sum * std::pow(-1.0/h,pow);\newline $\mathrm{\}}$ \\ \hline 
\end{tabular}

\item  Approximate calculation of the "pow" order differential of function

\begin{tabular}{|p{4.4in}|} \hline 
double N\_derivative(std::function$\mathrm{<}$double(double)$\mathrm{>}$ f, double pow, double point, unsigned long N) noexcept\newline $\mathrm{\{}$\newline     double sum = 0;\newline     for (size\_t k = 0; k $\mathrm{<}$ N; k++)//gofor(k, range(0, N))\newline     $\mathrm{\{}$\newline         sum += binomiald(pow, k) / std::exp(std::lgamma(k - pow + 1)) * std::pow(point, k - pow) * GL\_righthand\_derivative(f, k, point, N);\newline     $\mathrm{\}}$\newline     return sum;\newline $\mathrm{\}}$ \\ \hline 
\end{tabular}

\end{enumerate}

\noindent \textbf{Example 4.}

For example, investigate the transformation$D^{1/3} (x\wedge 2)$. The exact answer will be$\frac{\Gamma(3)}{\Gamma (8/3)} x^{5/3} $

 Table 1 - Comparison of the accuracy of fractional derivative calculation algorithms

\begin{tabular}{|p{0.4in}|p{1.3in}|p{1.3in}|p{1.3in}|} \hline 
Point & Exact answer & GL with 9 iteration & Eq \eqref{1.13} with 3 iteration \\ \hline 
0.2 & 0.0909213773601394 & 0.0909397004356575 & 0.0909213773601394 \\ \hline 
0.4 & 0.288657380135766 & 0.2887155522744711 & 0.288657380135766 \\ \hline 
0.6 & 0.5673722586492711 & 0.5674865992481914 & 0.5673722586492709 \\ \hline 
0.8 & 0.9164300577717989 & 0.9166147428001502 & 0.9164300577717989 \\ \hline 
1.0 & 1.329278600918967 & 1.329546485908133 & 1.329278600918967 \\ \hline 
1.2 & 1.801294640474853 & 1.801657649248871 & 1.801294640474852 \\ \hline 
1.4 & 2.328963807359781 & 2.329433155503858 & 2.328963807359781 \\ \hline 
1.6 & 2.909484075524463 & 2.910070413941038 & 2.909484075524462 \\ \hline 
1.8 & 3.540545570594309 & 3.541259084684707 & 3.540545570594309 \\ \hline 
\end{tabular}

Approximation of the solution of a specific differential equation.
\begin{equation} \label{1.18} 
D^{\alpha } u(x)+bx^{k-\alpha } \frac{d^{k} }{dx^{k} } u(x)=g(x) 
\end{equation} 

Equation \eqref{1.18} can be solved by using n first members of the sequence \eqref{1.13}
\begin{equation} \label{1.19} 
\sum _{k=0}^{n}\left(\begin{array}{c} {\alpha } \\ {k} \end{array}\right)\frac{x^{k-\alpha } }{\Gamma (k-\alpha +1)} \frac{d^{k} }{dx^{k} } u(x) +bx^{k-\alpha } \frac{d^{k} }{dx^{k} } u(x)=g(x) 
\end{equation} 

The resulting differential equation can be solved using standard methods. The resulting answer will contain an approximate solution and a set of functions that depend on constants. The solution will be closest when the constants are defined as zeros.
\newline
\noindent \textbf{Example 5.}

Consider the differential equation.
\[D^{0.5} u(x)+x^{0.5} \frac{d}{dx} u(x)=\frac{\ln (4x)}{\sqrt{\pi x} } +\frac{1}{\sqrt{x} } \] 

First, replace the fractional derivative with the first 3 elements of the sum of the equation \eqref{1.13}.
\[\sum _{k=0}^{2}\left(\begin{array}{c} {0.5} \\ {k} \end{array}\right)\frac{x^{k-0.5} }{\Gamma (k+0.5)} \frac{d^{k} }{dx^{k} } u(x) +x^{0.5} \frac{d}{dx} u(x)=\frac{\ln (4x)}{\sqrt{\pi x} } +\frac{1}{\sqrt{x} } \] 

Solve the equation with the help of software, for example, wolframalpha.
\[u(x)=0.170874+c_{2} x^{17.9686} +\frac{c_{1} }{x^{0.333915} } +\ln (x)\] 

The exact answer will be $u(x)=\ln (x)$ which means that we have gotten a solution that differs only by a constant when zeroing the constants in the answer.

\section{Discussion}

This work creates a number of issues that require future consideration. One of the interesting questions arises because of the high similarity of equality \eqref{1.13} in the form of $D^{\alpha } f(x)=\sum _{k=0}^{\infty }\left(\begin{array}{c} {\alpha } \\ {k} \end{array}\right)I^{k-\alpha }(1)D^{k} f(x) $ with the Newton binomial formula. Other formats for applying equality \eqref{1.13} can also be considered [8].
\newpage
\section*{References}

\begin{enumerate}
\item Zenyuk D.A., Orlov Yu.N. About application fractional Riemann-Liouville calculus for describing probability distributions // IPM preprints im. M.V. Keldysh. 2014. No. 18. 21 p. URL: http://library.keldysh.ru/preprint.asp?id=2014-18

\item Edmundo Capelas de Oliveira, Jos\'{e} Ant\'{o}nio Tenreiro Machado, "A Review of Definitions for Fractional Derivatives and Integral", Mathematical Problems in Engineering, vol. 2014, Article ID 238459, 6 pages, 2014. https://doi.org/10.1155/2014/238459

\item Chii-Huei Yu, " Study of Fractional Analytic Functions and Local Fractional Calculus, International Journal of Scientific Research in Science, Engineering and Technology(IJSRSET), Print ISSN : 2395-1990, Online ISSN : 2394-4099, Volume 8, Issue 5, pp.39-46, September-October-2021. Available at doi : https://doi.org/10.32628/IJSRSET218482 ~~

\item Madan, S. Some G P\'{o}lya gems from complex analysis. Reson 19, 323--337 (2014). https://doi.org/10.1007/s12045-014-0038-6

\item Puchta, Jan-Christoph. "On Fabry's gap theorem." Archivum Mathematicum 038.4 (2002): 307-309. $\mathrm{<}$http://eudml.org/doc/248932$\mathrm{>}$.

\item Osgood, Brad \& Wu, William. (2008). Falling Factorials, Generating Functions, and Conjoint Ranking Tables. 12.

\item Persechino, Andr\'{e}. (2020). An introduction to fractional calculus Numerical methods and application to HF dielectric response. Smart Materials and Structures. 9. 20. 10.7716/aem.v9i1.1192.

\item Ghorbani, H., Mahmoudi, Y., Saei, F.D. (2020). Numerical study of fractional Mathieu differential equation using radial basis functions. Mathematical Modelling of Engineering Problems, Vol. 7, No. 4, pp. 568-576. https://doi.org/10.18280/mmep.070409
\end{enumerate}

\end{document}